\documentstyle[11pt]{article}
\def\dfrac{\displaystyle\frac}
\def\f{\displaystyle}
\title{$SO(3)$ invariants of Seifert manifolds \\and their algebraic integrality
 \thanks{  Supported partially by the National Science
 Foundation of P. R. China, }}
\author{ Bang-He Li }
\date{}
\begin{document}
\maketitle
\baselineskip 18pt
\begin{abstract}

For Seifert manifold $M=X({p_1}/_{\f{q_1}},{p_2}/_{\f{q_2}}, \cdots,{p_n}/_
{\f{q_n}}),\;\tau^{'}_r(M)$ is calculated for all $r$ odd $\geq 3$. If $r$
is coprime to at least $n-2$ of $p_k$ (e.g. when $M$ is the Poincare
homology sphere), it is
proved that $(\sqrt {\dfrac{4}{r}}\sin \dfrac{\pi}{r})^{\nu}\tau^{'}_r(M)$ is
an algebraic integer in the r-th cyclotomic field, where $\nu$ is the first
Betti number of $M$. For the torus bundle obtained from trefoil knot with
framing $0$, i.e. $X_{\mbox{tref}}(0)=X(-2/_{\f{1}},3/_{\f{1}},6/_{\f{1}}),\; \tau^{'}_r$ is obtained in
a simple form if $3\mid\!\llap /r$, which shows in some sense that it is
impossible to generalize Ohtsuki's invariant to 3-manifolds being not
rational homology spheres.

\end{abstract}
\vskip 36pt
\section {Introduction}
\hskip\parindent
Consider Seifert manifolds of the form $M=
X({p_1}/_{\f{q_1}},{p_2}/_{\f{q_2}}, \cdots,{p_n}/_{\f{q_n}})$, where $p_k$
and $q_k$ are
not zero and coprime (the case of some $p_k$ or $q_k$ being zero is trivial,
 as explained in [GF]), while $\sum_k \dfrac{q_k}{p_k}=0$ (i.e. $M$ is not
 a rational homology sphere) is allowed.

In the case of $\sum_k \dfrac{q_k}{p_k}\not=0$, $r$ being a prime
and $r\mid\!\llap / {p_k},\,
r\mid\!\llap / {q_k},\, k=1,2,\cdots,n$, Rozansky [R1] has obtained a simpler
formula for $\tau^{'}_r(M)$ defined by Kirby and Melvin [KM1].
Except obtaining a general formula for all odd $r\geq 3$, we obtain
simpler formula in the case of $r$ being coprime to all $p_k$. And after
changing Rozansky's formula, we see that it becomes our form
in his case.

Now, let us introduce some notations. $P=\prod_k p_k,\, H=P\sum_k
\dfrac{q_k}{p_k},\,c_k=(r, p_k)$ being the commom factor, $s(q,p)$ is the
Dedekind sum, and $(\dfrac{a}{b})$ is the Jacobi symbol for odd $b>0$, while
the ratio of $c$ to $d$ will be written as $\frac{c}{d}$ or $c/_{\f{d}}\,$.
$p^*_k$and $q^*_k$ are any pair satisfying
$$p^*_kp_k+q^*_kq_k=1$$
while for $c_k=(r,p_k),\, ({p_k}/_{\f{c_k}})'$ and $({r}/_{\f{c_k}})'$ are any pair defined by
$$({p_k}/_{\f{c_k}})'\frac{p_k}{c_k}+({r}/_{\f{c_k}})'\frac{r}{c_k}=1$$
The functions $\varepsilon(r)$ is defined by
 $\varepsilon(r)=1$, if $r\equiv 1\pmod 4$,
 $\varepsilon(r)=i=\sqrt {-1}$, if $r\equiv -1 \pmod 4$.

And if $c_k=1, 12s^{\surd}(q_k,p_k)\equiv (12p_k(q_k,p_k))p_k' \pmod r$.
Moreover, $\chi^{k,\pm}(j)$ is defined by
$$\chi^{k,\pm}(j)=\cases{\pm 1, & if $\,j\mp q^*_k\equiv 0 \pmod {c_k}$\cr
0,&$\;$ otherwise\cr}$$
and $e_a$ denotes $e^{\frac{2\pi i}{a}}$, $A=e_r^{\frac{1\mp r}{4}}$ for
$r\equiv \pm 1\pmod 4$. We have
\smallskip

{\bf Theorem 1}. With the notations as above, and assume all $q_k>0, r$ odd $\geq 3$,
then
$$\begin{array}{ll}
\tau^{'}_r(X({p_1}/_{\f{q_1}},&{p_2}/_{\f{q_2}}, \cdots,{p_n}/_{\f{q_n}}))
=(\sqrt{1\over r}{\sin{\pi\over r}})^{1-sign |H|} (-2\varepsilon(r)\sqrt r)
^{-sign |H|}\times\\\\
&(sign P (-sign \dfrac{H}{P}+1-sign |H|))
^{\frac{r+1}{2}}(A^2-A^{-2})^{-2+sign |H|}\times\\\\
& A^{\displaystyle -3sign \dfrac{H}{P}+\sum^n_{k=1}(-12 s(q_k,p_k)+
\dfrac{q_k+q^*_k}{p_k}-({p_k}/_{\f{c_k}})' {{q^*_k}\over{c_k}}
-(r/_{\f{c_k}})'{r\over{c_k}}p^*_k q^*_k)}\times\\\\
&\prod^n_{k=1}{\displaystyle ((-1)^{\dfrac{r-1}{2}\dfrac{c_k-1}{2}}\sqrt{c_k}
\varepsilon(c_k)(\dfrac{{p_k}/_{\f{c_k}}}{r/_{\f{c_k}}})(\dfrac{q_k}{c_k}))}\times\\\\
&\sum^r_{j=1}(A^{2j}
-A^{-2j})^{2-n}A^{-(\sum^n_{k=1}{\displaystyle ({p_k}/_{\f{c_k}})'
\dfrac{q_k}{c_k})j^2}}\times\\\\
&\prod ^n_{k=1}\sum_{\pm}\chi^{k,\pm}(j)
A^{\displaystyle \pm 2(\dfrac{1}{c_k}({p_k}/_{\f{c_k}})'+p^*_k
({r}/{c_k})'\dfrac{r}{c_k})j}\end{array}$$
\vskip 36pt
{\bf Corollary}. If $r$ is coprime to all $p_k$, i.e. $c_k=1$ for
$k=1,\cdots,n$, then
$$\begin{array}{rl}
\tau^{'}_r(X({p_1}/{q_1},{p_2}/{q_2},&\cdots,{p_n}/{q_n}))
=(\sqrt{1\over r}\sin{\pi\over r})^{1-sign |H|} (-2\varepsilon(r)\sqrt r)
^{-sign |H|}\times\\\\
& (-sign \dfrac{H}{P}+1-sign |H|)
^{\frac{r+1}{2}}(A^2-A^{-2})^{-2+sign |H|}(\dfrac{|P|}{r})sign P\times\\\\
&A^{-3sign \dfrac{H}{P}+P'H-12\sum^n_{k=1} s^{\surd}(q_k,p_k)}\times\\\\
&\sum^r_{j=1}(A^{2j}-A^{-2j})^{2-n} A^{-P'Hj^{2}}\prod^n_{k=1}
(A^{2p_k^{'}j}-A^{-2p_k^{'}j})
\end{array}$$
where $P'P\equiv 1\pmod r$
\smallskip

{\bf Remark 1}. In Theorem 1, the assumption of all $q_k >0$ is necessary,
otherwise the formula would be changed to a more complicated form. It can be seen from the proof.
  While for the corollary, we need not this assumption.
\smallskip

{\bf Theorem 2}. If $r$ is coprime to ar least $n-2$ of $p_k$, then
$$\tau^{'}_r (X({p_1}/{q_1},{p_2}/{q_2}, \cdots,{p_n}/{q_n})
(\sqrt{1\over r}\sin{\pi\over r})^{sign |H|-1}$$ is an algebraic  interger
in the $r$-th cyclotomic field.
\smallskip

{\bf Remark 2}. We actually prove that $\Theta_r(M)$ and $\xi_{r}(M)$
defined in [BHMV] and [Li]
respectively are all algebraic integers in this case.
\smallskip

From Theorem 2, we see that for all rational homology sphere of the form\\
$X(2/{q_1},\; {p_2}/{q_2},{p_3}/{q_3}),\;\tau_r^{'}$
is an algebraic integer for all odd $r\geq 3$. Especially it is true for all RHS obtained by
Dehn surgery on $(2,n)$ torus knot, including those considered by R. Lawsence in
[La]. Interestingly the Poincare homology sphere $\sum (2,3,5)$ and all the
Brieskon homology spheres of the form $\sum (2^n, p,q)$ are among them.

For the torus bundle over $S^1$ obtained by the monodromy matrix
$$\left(\begin{array}{cc}
0& -1\\
1&1\end{array}\right) \in SL(2,Z)$$
which is the  Seifert manifold $X(-2/_{\f{1}}, 3/_{\f{1}}, 6/_{\f{1}})$ or
$X_{\mbox{tref}}(0)$, i.e. gotten by doing surgery on left-handed trefoil
knot with framing  $0$, we calculate $\tau^{'}_r$ further and get
$$\tau^{'}_r(X_{\mbox{tref}}(0))=\cases{0,& if $r\equiv 1\;\;\,\pmod 3$\cr
\dfrac{-\sqrt r}{2sin{\pi\over r}}e^{-\frac{2\pi}{r}i},&
if $r\equiv -1\pmod 3$\cr}$$
For rational homology 3-sphere $M$, there is Ohtsuki's invariant
$$\tau (M)=\sum^{\infty}_{n=0}\lambda _n(M)(t-1)^n \in Q[[t-1]]$$
where $\lambda_n(M)$ is determined by $\tau^{'}_r(M)$ with sufficiently large
prime $r$. By the Conjecture of R. Lawsence [La] proved recently by Rozansky
[R2], if $|H_1(M,Z)|\not= 0\pmod r$
then the cyclotomic series $\sum^{\infty}_{n=0}\lambda^r_n(M)h^n$ converges
r-adicly to
$$(\dfrac{|H_1(M,Z)|}{r}) |H_1(M,Z)|\tau^{'}_r(M)$$
where $h=e^{{2\pi i}\over r}-1$.
Therefore for any suitably large prime $r$, $\tau^{'}_r(M)$ is determined by
 any sequence of $\{\tau^{'}_{r_n}(M)\}$, where $r_n$ is prime and
 $\lim_{n\rightarrow\infty}r_n=\infty$. While for $M=X_{\mbox{tref}}(0)$,
 we see from the formula above that this property does not hold. Since we
 have such a sequence of $\{\tau^{'}_{r_n}(M)\}$ with $\tau^{'}_{r_n}(M)=0$,
but there exists abitrarily large prime $r$ with $\tau^{'}_{r}(M)\not=0$.
It thus seems that generalizing Ohtsuki's invariant to general 3-manifolds
is impossible.
\bigskip

{\bf Remark 3}. For lens spaces, $\xi_r$ which is equivalent to $\tau_r^{'}$
 has been obtained in [LL1], and explicit formula for $\tau_r^{'}$ is obtained
 in [LL2] in order to calculate Ohtsuki's invariant. An example is given in
 [LL3] to show that Ohtsuki's invariant does not determine all $\tau_r^{'}$.

\section { Proof of Theorem 1 and its corollary}

{\bf 2.1  Reducing to $\xi_{r}(M,e_{r})$}

\hskip\parindent
Recall from [Li] that
$$\tau_{r}^{'}(M)=(\sqrt{4\over {r}}\sin{ \pi\over r})^{\nu}\Theta_{r}(M,\pm ie_{4r}),~~~~
~~3\leq r\equiv\pm 1 \pmod 4          $$
where $e_{a}$ stands for $e^{2\pi i\over a}, \;\nu$ is the first Betti number
of $M$ , and $\Theta_{r}$ is defined in [BHMV]. And
$$\Theta_{r}(M,\pm ie_{4r})=2^{-\nu}\xi_{r}(M,\mp ie_{4r})=2^{-\nu}\xi_{r}~
(M,e_{r}^{1\mp r\over 4}),~~~~~3\leq r\equiv\pm 1 \pmod 4  $$
where $\xi_{r}$ is defined in [Li] and the formula above is also proved
there. Therefore
$$\tau_{r}^{'}(M)=(\sqrt{1\over r}\sin{ \pi\over r})^{\nu}\xi_{r}~
(M,e_{r}^{1\mp r\over 4})    $$

Our working line is then to calculate $\xi_{r}(M,e_{r})$ first , then use the
Galois automorphism of the r-th cyclotomic field sending $e_{r}\rightarrow~
e_{r}^{1\mp r\over 4} $ to  get $\xi_{r}(M,e_{r}^{1\mp r\over 4})  $,
hence $\tau_{r}^{'}(M) $.

For a rational number $\alpha$, denote by
$<m_{l},m_{l-1},\ldots, m_{1}>$ its continued fraction expression
$$\alpha=m_l-{{1}\over\displaystyle m_{l-1}-
{\strut\ \atop{\ddots -\displaystyle{\strut{1}
\over\displaystyle m_2-{\strut{1}\over\displaystyle {m_1}}}}}}$$

let $p_k/_{\f{q_k}}= < m_{k,l_k,} m_{k,l_{k-1},}\cdots, m_{k,l_1} > $, then
$X(p_1/_{\f{q_1}},\cdots, p_n/_{\f{q_n}})$ is represented by the framed link:

\bigskip

where a dot $\bullet$ lebbelled with a number $m$ means an unknot with
framing $m$, and two dots connected by a line means that 2 relevant
unknots form a Hopf link.

This shows that $X(p_1/_{\f{q_1}},\cdots, p_n/_{\f{q_n}})$ is a plumed manifold. By using the formula
 in [LL4] for $\Theta_r$ of plamed manifolds together with relation between
 $\Theta_r$ and $\xi_r$ in [Li], we get
$$\xi_r(X(p_1/q_1,\cdots, p_n/q_n), e_r)=I/II$$
where
$$I=(e_r^2-e_r^{-2})^{-(N+1)}e_r^{-\sum_{k,l}m_{k,l}}\sum^r_{j=1}
(e_r^{2j}-e_r^{-2j})^{2-n}\prod^n_{k=1}S_{k,l_k}(j)$$
with $N$ being the number of components of
the framed link, i.e. $N=\sum_k l_k+1$, and
$$S_{k,l_k}(j)=\sum^r_{j_1,\cdots,j_{l_k}=1}
e_r^{\sum^{l_k}_{l=1}m_{k,l}\,j^2_l}
(e_r^{2j_1}-e_r^{-2j_1})(e_r^{2j_{l_k}j}-e_r^{-2j_{l_k}j})\prod_{t=1}
^{l_k-1}(e_r^{2j_tj_{t+1}}-e_r^{-2j_tj_{t+1}})$$
and
$$II=s_+^{b^+}s_-^{b^-},\;\;\; s_+=\dfrac{-2e_r^{-3}}{e_r^2-e_r^{-2}}\sqrt r
\varepsilon(r),\;\;\;s_-=\bar{s}_+$$
with $b_+$ and $b_-$ the numbers of positive and negative eigen values of the
linking matrix respectively, and $\bar{s}_+$ means the complex conjugete of
$s_+$.
\bigskip

{\bf 2.2. Good continued franction expression.}

\hskip\parindent
We call $<m_{l},m_{l-1},\ldots, m_{1}>$ a good continued framction expression
of $\alpha$, or simply a good expression of $\alpha$, if $m_l=[\alpha]+1$
with $[\alpha]$ the integer part of $\alpha$, and
$<m_{l-1},\ldots, m_{1}>=([\alpha]+1-\alpha)^{-1}$ such that
 $m_{1},m_{2},\ldots, m_{l-1}$ are all $\geq 2\,$ if $\,[\alpha]\not=\alpha$,
and $<m_{l-1},\ldots, m_{1}>=<1>$ if $[\alpha]=\alpha$.

In [LL1], for $\alpha=p/q$ with $p>0,\;q>0$, and $\alpha=< m_{l},m_{l-1},
\ldots,m_{1}>$ with all $m_j\geq 2,\; N_{j,i}$ for $1\leq i\leq j\leq l$ is
defined to be the numerator of $<m_{j},\ldots,
m_{i}>$, then $$N_{l,1}=p,\;\;N_{l-1,1}=q$$

Let $p_k/_{\f{q_k}}=<m_{k,l_k},\cdots,m_{k,1}>$ be a good expression, then
for $ i\leq j\leq l_k-1$, since $<m_{k,j},\cdots,m_{k,i}>$ is positive, we can still define
$N_{k;j,i}$ as its numerator, and in the assumption of $q_k>0$, we have
$N_{k;l_k-1,1}=q_k$. Thus althrough $p_k$ can be negative, we still have $p_k=
N_{k;l_k,1}$ being the numerator of $<m_{k,l_k},\cdots,m_{k,1}>$ and all
results concerning $N_{j,i}$ are still true.
\bigskip

{\bf 2.3  Calculation for $S_{k,l_k}(j)$}

\hskip\parindent
By Lemma 4.12, 4.20 and some other lemmas in \S 4 of [LL1], we have
$$\begin{array}{rl}S_{k,l_k}(j)=&(-2\sqrt r\varepsilon(r))^{l_k}\sqrt{c_k}
\varepsilon(c_k)(\dfrac{p_k/_{\f{c_k}}}{r/_{\f{c_k}}})(\dfrac{q_k}{c_k})\times\\
&(-1)^{\frac{r-1}{2}\frac{c_k-1}{2}}\sum_{\pm}\chi^{k,\pm}
(j)e_r^{-({p_k}/_{\f{c_k}})'\dfrac{q_k}{c_k}(j\mp q^*_k)^2-p^*_k(q^*_k\mp 2j)}
\end{array}$$
where $q^*_k=N_{k;l_k,2}$ and $p^*_k=-N_{k;l_k-1,2}$ is a spicial choice of
$q^*_k$ and $p^*_k$ satisfying $q^*_kq_k+p^*_kp_k=1$.

Now
$$-({p_k}/_{\f{c_k}})'\dfrac{q_k}{c_k}(j\mp q^*_k)^2-p^*_k(q^*_k\mp 2j)$$
$$=-({p_k}/_{\f{c_k}})'\dfrac{q_k}{c_k}j^2\pm 2(({p_k}/_{\f{c_k}})'\dfrac{q_k q^*_k}{c_k}
+p^*_k)j-({p_k}/_{\f{c_k}})'\dfrac{q_k}{c_k}q^*_kq^*_k-p^*_kq^*_k$$
and
$$\begin{array}{rl}({p_k}/_{\f{c_k}})'\dfrac{q_k q^*_k}{c_k}+p^*_k=&
({p_k}/_{\f{c_k}})'\dfrac{1-p_k p^*_k}{c_k}+p^*_k=
\dfrac{1}{c_k}({p_k}/_{\f{c_k}})'-p^*_k(1-(r/_{\f{c_k}})'\dfrac{r}{c_k})+p^*_k\\=&
\dfrac{1}{c_k}({p_k}/_{\f{c_k}})'+p^*_k(r/_{\f{c_k}})'\dfrac{r}{c_k},\end{array}$$
$$\begin{array}{rl}-({p_k}/_{\f{c_k}})'\dfrac{q_k}{c_k}q^*_kq^*_k-p^*_kq^*_k=&
-({p_k}/_{\f{c_k}})'\dfrac{q^*_k}{c_k}(1-p^*_kp_k)-p^*_kq^*_k\\=&
-({p_k}/_{\f{c_k}})'\dfrac{q^*_k}{c_k}+p^*_kq^*_k(1-(r/_{\f{c_k}})'
\dfrac{r}{c_k})-p^*_kq^*_k\\
=&-({p_k}/_{\f{c_k}})'\dfrac{q^*_k}{c_k}-p^*_kq^*_k(r/_{\f{c_k}})'
\dfrac{r}{c_k}\end{array}$$
Thus
$$\begin{array}{rl}S_{k,l_k}(j)=&(-2\sqrt r\varepsilon(r))^{l_k}\sqrt{c_k}
\varepsilon(c_k)(\dfrac{p_k/_{\f{c_k}}}{
r/_{\f{c_k}}})(\dfrac{q_k}{c_k})(-1)^{\dfrac{r-1}{2}\dfrac{c_k-1}{2}}
e_r^{-({p_k}/_{\f{c_k}})'\dfrac{q^*_k}{c_k}-p^*_kq^*_k(r/_{\f{c_k}})'\dfrac{r}{c_k}}\times\\
&e_r^{-({p_k}/_{\f{c_k}})'\dfrac{q_k}{c_k}j^2}
\sum_{\pm}\chi^{k,\pm}(j)
e_r^{\pm 2(\dfrac{1}{c_k}({p_k}/_{\f{c_k}})'+p^*_k(r/_{\f{c_k}})'\dfrac{r}{c_k})j}
\end{array}$$
\bigskip

{\bf 2.4. Calculation for $\xi_r(X(p_1/_{\f{q_1}},\cdots, p_n/_{\f{q_n}}), e_r)$}

\hskip\parindent
First by a diagonalization procedure for the quatratic form of the linnk matrix,
we see that
$$b_++b_-=\cases {N, & if $\sum \dfrac{q_k}{p_k}\not=0$\cr
N-1, & if $\sum \dfrac{q_k}{p_k}=0 $\cr}$$
$$b_-=\mbox{the number of the negative elements in the set}\;
\{p_1/q_1,\cdots, p_n/q_n, -\sum \dfrac{q_k}{p_k}\}, $$
since we use a food expression for every $p_k/q_k$. Thus by the assumption
$q_k>0$, we have
$$\begin{array}{rl}b_++b_-&=N-1+sign |H|\cr
b_-&=\mbox{the number of the negative elements in the set}\;
\{p_1,\cdots, p_n, -\dfrac{H}{P}\}
\end{array}$$
and $(-1)^{\frac{r+1}{2}b_-}=(sign P)^{\frac{r+1}{2}}(-sign
\dfrac{H}{P}+1-sign |H|)^{\frac{r+1}{2}}$.

Now
$$s_+^{-b_+}s_-^{-b_-}=(e_r^2-e_r^{-2})^{b_++b_-}(-2\sqrt r)^{-b_+-b_-}
(-1)^{b_-}\varepsilon(r)^{-b_++b_-}$$
Thus
$$\xi_r(X(p_1/_{\f{q_1}},\cdots, p_n/_{\f{q_n}}), e_r)=I/II=G_1G_2$$
where
$$\begin{array}{ll}G_1=&e_r^{-\sum_k ({p_k}/_{\f{c_k}})'\frac{q^*_k}{c_k}-
\sum _{k}p^*_kq^*_k(r/_{\f{c_k}})'\dfrac{r}{c_k}}
\sum^r_{j=1}(e_r^{2j}-e_r^{-2j})^{2-n} e_r^{-\sum_k ({p_k}/_{\f{c_k}})'
\dfrac{q_k}{c_k}j^2}\times\\\\&
\prod_{k}\sum_{\pm}\chi^{k,\pm}(j)
e_r^{\pm 2(\dfrac{1}{c_k}({p_k}/_{\f{c_k}})'+p^*_k(r/_{\f{c_k}})'\dfrac{r}{c_k})j}
\end{array}$$
and
$$\begin{array}{ll}G_2=&e_r^{3(b_+-b_-)-\sum_{k,l}m_{k,l}}(e_r^2-e_r^{-2})^{b_++b_--N-1}
(-2\sqrt r)^{-b_+-b_-}(-1)^{b_-}\times\\\\
&\varepsilon(r)^{b_{-}-b_+}(-2{\sqrt r}\varepsilon(r))^{\sum_k l_k}
\prod_k(\sqrt{c_k}(\dfrac{p_k/_{\f{c_k}}}{
r/_{\f{c_k}}})(\dfrac{q_k}{c_k})(-1)^{\frac{r-1}{2}\frac{c_k-1}{2}}\varepsilon(c_k))
\end{array}$$

Notice that
$$b_+-b_-=-sign\dfrac{H}{P}+\sum_k(l_k-1+sign p_k)$$
and
$$3(l_k-1+sign p_k)-\sum^{l_k}_{l=1}m_{k,l}
=-12s(q_k,p_k)+\dfrac{q_k+q^*_k}{p_k}$$
(cf.[KM2]), where $q^*_k=N_{k;l_k,2}$. Also we have $\sum_k l_k=N-1$ and
$$(-1)^{b_-}\varepsilon(r)^{-(b_+-b_-)+\sum_k l_k}=
(-1)^{b_-}\varepsilon(r)^{2b_--sign |H|}=
(-1)^{\frac{r+1}{2}b_-}\varepsilon(r)^{-sign |H|}
$$
So,
$$\begin{array}{ll}G_2=&e_r^{-3sign \frac{H}{P}+
\sum_k(-12s(q_k,p_k)+\frac{q_k+q^*_k}{p_k})}
(e_r^2-e_r^{-2})^{-2+sign |H|}(-2\sqrt r\varepsilon(r))^{-sign |H|}\times\\\\&
(sign P)^{\frac{r+1}{2}}(-sign \frac{H}{P}+1-sign |H|)^{\frac{r+1}{2}}
\prod_k(\sqrt{c_k}(\dfrac{p_k/_{\f{c_k}}}{r/_{\f{c_k}}})(\dfrac{q_k}{c_k})
(-1)^{\frac{r-1}{2}\frac{c_k-1}{2}}\varepsilon(c_k))
\end{array}$$
and $\xi_r(M,e_r)$ is obtained.
\bigskip

{\bf 2.5. The Galois automorphism}

\hskip\parindent
Consider the Galois automorphism sending $e_r$ to $A=e_r^{\frac{1\mp r}{4}}$
for $r\equiv \pm 1 \pmod 4$. In the formula for $\xi_r(M,e_r)$, since
$$G_1=\sum^r_{j=1}(e_r^{2j}-e_r^{-2j})^{2-n}\prod_k\sum_{\pm}\chi^{k,\pm}(j)
e_r^{-(p_k/_{\f{c_k}})'\dfrac{q_k}{c_k}(j\mp q^*_k)^2-p^*_k(q^*_k\mp 2j)}$$
and when $\chi^{k,\pm}(j)\not=0,\;c_k|\,j\mp q^*_k$, we see that the image of
$G_1$ under the automorphism is obtained via replacing $e_r$ by $A$.

For any positive factor $c$ of $r$, we have
$$\sqrt c \varepsilon(c)=\sum^c_{x=1}e_c^{x^2}=\sum^c_{x=1}
e_r^{\frac{r}{c}x^2}$$
So the image of $\sqrt c\varepsilon (c)$ is
$$\sum^c_{x=1}e_r^{\frac{1\mp r}{4}\frac{r}{c}x^2}=\sum^c_{x=1}
e_c^{\frac{1\mp r}{4}x^2}=\sqrt c \varepsilon(c)(\frac{{(1\mp r)}/4}{c})$$
$4\cdot \frac{1\mp r}{4}\equiv 1\pmod r$ implies
$4\cdot \frac{1\mp r}{4}\equiv 1\pmod c$.
This means that $ \frac{1\mp r}{4}\pmod c$ is a square.
Hence
$$\Bigg{(}\frac{ (1\mp r)/{4}}{c}\Bigg{)}=1$$
and the image of $\sqrt c \varepsilon(c)$ is $\sqrt c \varepsilon(c)$. Thus
the image of $G_2$ is obtained by putting $A$ in the place of $e_r$.
\bigskip

{\bf 2.6 Independentness of the choice of $q^*_k$ and $p^*_k$}

\hskip\parindent
So far in the formula for $\xi_r(M,e_r)$, $q^*_k$ and $p^*_k$ depends on the
good exrpression of $p_k/q_k$. If $q^*_k$ is changed, it must become
$q^*_k+m p_k$ for some $m\in Z$, and then $p^*_k$ becomes $p^*_k-mq_k$.
Look at the change of $G_1$, we see that
$\chi^{k,\pm}(j)$ does not change.\\
While $
-({p_k}/_{\f{c_k}})'\dfrac{q_k}{c_k}(j\mp q^*_k)^2-p^*_k(q^*_k\mp 2j)$ changes to
$$\begin{array}{ll}
&-({p_k}/_{\f{c_k}})'\dfrac{q_k}{c_k}(j\mp q^*_k\mp mp_k)^2
-(p^*_k-mq_k)(q^*_k\mp 2j+mp_k)\\
=&-({p_k}/_{\f{c_k}})'\dfrac{q_k}{c_k}(j\mp q^*_k)^2-p^*_k(q^*_k\mp 2j)+X+Y
\end{array}$$
where
$$\begin{array}{ll}
X&=-m({p_k}/_{\f{c_k}})'\dfrac{p_k}{c_k}(mp_k+2(q^*_k\mp j))q_k\\\\
Y&=-mp^*_kp_k-mq^*_kq_k+m^2p_kq_k+2mq_k(q^*_k\mp 2j)
\end{array}$$
Since $({p_k}/_{\f{c_k}})'\dfrac{p_k}{c_k}\equiv 1\pmod {\frac{c_k}{r}}$, and
$mp_k+2(q^*_k\mp j)\equiv 0\pmod r$,
$$X\equiv -m(mp_k+2(q^*_k\mp j))q_k \pmod r$$
Thus
$$X+Y\equiv -mp^*_kp_k-mq^*_kq_k\equiv -m\pmod r$$
and $G_1$ changes to $e_r^{-m}G_1$.

Now it is obvious that $G_2$ changes to $e^m_rG_2$. So $\xi_r(M,e_r)$ does not depends on
the special choice of $q^*_k$ and $p^*_k$, and the proof for Theorem 1 is complete.
\bigskip

{\bf 2.7. Proof of the Corollary}

\hskip\parindent
Since
$$-12s(q_k,p_k)+\dfrac{q_k+q^*_k}{p_k}=3(l_k-1+sign p_k)
-\sum^{l_k}_{l=1}m_{k,l}$$
and
$$-12s^\surd(q_k,p_k)\equiv 3(l_k-1+sign p_k)-\sum^{l_k}_{l=1}m_{k,l}
-p^{'}_kq^*_k-p^{'}_kq_k\;\pmod r$$
$$\sum_k p^{'}_k q_k=P'H$$
$$\prod_k (\frac{p_k}{r})=(\frac{P}{r})=(\frac{|P|}{r})(sign P)^{\frac{r-1}{2}}$$
we are done.
\bigskip
\section {Comparision with the formula of Rozansky}
\hskip\parindent
Under the assumption of $p_k, \,q_k\not\equiv 0\pmod r,\, H\not=0$, and $r$
 being prime, Rozansky's formula for $\tau^{'}_r(X
 ({p_1}/_{\f{q_1}},{p_2}/_{\f{q_2}}, \cdots,{p_n}/_{\f{q_n}}))$ is
$$\begin{array}{ll}\tau^{'}_r=\frac{i}{2\sqrt r}&e
^{\frac{i\pi}{4}sign \frac{H}{P}(\varepsilon(r)^2+3\frac{r-2}{r})}\\\\
&\times (\frac{|P|}{r})sign P e_r^{4'P'H-3\sum^n_{k=1}s^\surd (q_k,p_k)}
(e_r^{2'}-e_r^{-2'})^{-1}\\\\
&\times \sum_{{0\leq \beta<2r}, {\beta\in 2Z+1}}
e_r^{-4'P'H\beta^2}
(e_r^{2'\beta}-e_r^{-2'\beta})^{2-n}\prod_{k=1}^n (e_r^{2'p^{'}_k\beta}-
e_r^{-2'p^{'}_k\beta})\end{array}$$

Let $\beta=2(\alpha-2')+1$, then $\beta\equiv 2\alpha\pmod r$. And let $\alpha=2'j$, then we have
$$\begin{array}{ll}\sum_{{0\leq \beta<2r}, {\beta\in 2Z+1}}& e_r^{-4'P'H\beta^2}
(e_r^{2'\beta}-e_r^{-2'\beta})^{2-n}\prod_{k=1}^n (e_r^{2'p^{'}_k\beta}-
e_r^{-2'p^{'}_k\beta})\\
&=\sum^r_{\alpha=1}e_r^{-P'H\alpha^2}(e_r^{\alpha}-e_r^{-\alpha})
^{2-n}\prod _{k=1}^n(e_r^{p^{'}_k\alpha}-e_r^{-p^{'}_k\alpha})\\
&=\sum_{j=1}^r(A^{2j}-A^{-2j})^{2-n}A^{-P'Hj^2}\prod^n_{k=1}(A^{2p^{'}_kj}-
A^{-2p^{'}_kj})\end{array}$$

It can be checked that in his formula the term
$$ie^{\frac{i\pi}{4}sign \frac{H}{P}(\varepsilon(r)^2+3\frac{r-2}{r})}=
(sign \frac{H}{P})^{\frac{r+1}{2}}\varepsilon(r)A^{-3sign \frac{H}{P}}$$
And in our formula the term
$$\begin{array}{rl}
(-\varepsilon(r))^{-sign |H|}(-sign\frac{H}{P}+1-sign|H|)^{\frac{r+1}{2}}&=
-\varepsilon(r)(-1)^{\frac{r-1}{2}}(-sign\frac{H}{P})^{\frac{r+1}{2}}\\
&=\varepsilon(r)(sign\frac{H}{P})^{\frac{r+1}{2}}\end{array}$$
Hence two formulas coincide.
\bigskip
\section {Proof of Theorem 2}
\smallskip
\hskip\parindent
Theorem 2 is equivalent to algebraic integrality of $
\xi_r(X(p_1/_{\f{q_1}},\cdots, p_n/_{\f{q_n}}), e_r)$,
when $r$ is coprime to at least $n-2$ of $p_k$,

We assume $c_k=1$ for $3\leq k\leq n$. Then
$$\xi_r(M,e_r)=(2\varepsilon(r)\sqrt r)^{-sign |H|} (e_r^2-e_r^{-2})
^{-2+sign |H|}\varepsilon(c_1)\sqrt {c_1}\varepsilon(c_2)\sqrt {c_2}F_1F_2$$
where $F_1$ is an algebraic integer, and
$$F_2=\sum^r_{j=1}(K_{+}(j)+K_{-}(j))K(j)$$
with $$K_{\pm}(j)=\chi^{1,\pm}(j)e_r^{-(p_1/_{\f{c_1}})'\dfrac{q_1}{c_1}
(j\mp q^*_1)^2\pm 2p^{*}_1(j\mp q^*_1)}$$ and
$$\begin{array}{ll}
K(j)=&(\sum_{\pm}\chi^{2,\pm}(j)e_r^{-(p_2/_{\f{c_2}})'\dfrac{q_2}{c_2}
(j\mp q^*_2)^2\pm 2p^{*}_2(j\mp q^*_2)})\times\\\\
&e_r^{(-\sum^n_{k=3}p^{'}_kq_k)j^2}(e_r^{2j}-e_r^{-2j})^{2-n}
\prod^n_{k=3}(e_r^{2p^{'}_kj}-e_r^{-2p^{'}_kj})
\end{array}$$

Notice that $K(j)$ and $K_{\pm}(j)$ are all functions of $j$ with period $r$,
and
$$\chi^{k,\pm}(-j)=-\chi^{k,\mp}(j)$$
Hence
$$\begin{array}{ll}\sum^r_{j=1}K_{+}(j)K(j)&=\sum^r_{j=1}K_{+}(-j)K(-j)\\
&=\sum^r_{j=1}(-\chi^{1,-}(j)(\sum_{\pm}-\chi^{2,\mp}(j)
e_r^{-(p_2/_{\f{c_2}})'\dfrac{q_2}{c_2}
(j\pm q^*_2)^2\mp 2p^*_2(j\pm q^*_2)})\\\\&\times
e_r^{(-\sum^n_{k=3}p^{'}_kq_k)j^2}(e_r^{2j}-e_r^{-2j})^{2-n}
\prod^n_{k=3}(e_r^{2p^{'}_kj}-e_r^{-2p^{'}_kj}))\\\\
&=\sum^r_{j=1}K_{-}(j)K(j)\end{array}$$
and
$$F_2=2\sum^r_{j=1}K_{+}(j)K(j)$$

Let
$$\begin{array}{ll}F_2^{\pm}=&2\sum_{j=1}^r K_+(j)\chi^{2,\pm}(j)
e_r^{-(p_2/_{\f{c_2}})'\dfrac{q_2}{c_2}
(j\mp q^*_2)^2\pm 2p^*_2(j\mp q^*_2)}\\\\&\times
e_r^{(-\sum^n_{k=3}p^{'}_kq_k)j^2}(e_r^{2j}-e_r^{-2j})^{2-n}
\prod^n_{k=3}(e_r^{2p^{'}_kj}-e_r^{-2p^{'}_kj})\end{array}$$

Consider the set
$$S^{\pm}=\{j/\;\chi^{1,+}(j)\chi^{2,\pm}(j)\not=0\} \subset Z_r$$
If $S^{\eta}$ is empty, then $F_2^{\eta}=0$,
where $\eta=+$ or $-$. Now assume $S^{\eta}
\not=\emptyset$, and let $c$ be the least common multiple of $c_1$ and $c_2$. Then
$c|r$, and $S^{\eta}
$ is a residue class of $Z_r\pmod{Z_c}$, i.e.
$$S^{\eta}=\{\,a_\eta +xc / x\in Z_{r\over c}\}
$$
for some $a_\eta \in Z_r$.

Since
$$(e_r^{2j}-e_r^{-2j})^{2-n}\prod^n_{k=3}(e_r^{2p^{'}_kj}-e_r^{-2p^{'}_kj})
=\sum^m_{s=1}e_r^{2b_s j}$$
for some integers $b_1,\cdots, b_m$, we see that
substituting $j$ with $a_\eta +xc$ yields
$$
F_2^{\eta}=2\eta\sum_{s=1}^m(e_r^{\beta_{\eta, s}}
\sum_{x=1}^{r/_{\f{c}}}e_r^{c(Dx^2+E_{\eta,s}x)})
$$
where $\beta_{\eta,s}$, $E_{\eta,s}$ are integers depending on $\eta$ ans $s$,
while $D$ is an integer independent of $\eta$ ans $s$.
\bigskip

Let $d=(D, r/_{\f{c}})$, we have by Theorem 2.1 and 2.2 in [LL1]
$$\sum_{x=1}^{r/_{\f{c}}}e_{r/_{\f{c}}}^{Dx^2+E_{\eta,s}x}=\cases {r/_{\f{c}}, &if $D=0$,
and ${r\over c} |E_{\eta,s}$\cr
0, &if $D=0$, and ${r\over c}\mid\!\llap / E_{\eta,s}$ \cr&or $D\not=0$ and
$d \mid\!\llap / E_{\eta,s}$\cr
d\sum^{r/_{\f{cd}}}_{x=1}e_{r/_{\f{cd}}}^{\frac{D}{d}x^2+\frac{E_{\eta,s}}{d}x}=\pm d\varepsilon(r/_{\f{cd}})\sqrt
{r/_{\f{cd}}}\,e_r^{F_{\eta,s}}, & if $D\not=0$ and $d|E_{\eta,s}$\cr}$$
where $F_{\eta,s}$ is an integer. Thus
$$F_2^{\eta}=\cases {2{r\over c} c^{\eta}, &if $D=0$\cr
2d\varepsilon(r/_{\f{cd}})\sqrt{r/_{\f{cd}}}\,c^{\eta}, &if $D\not=0$\cr}$$
for some algebraic integer $c^{\eta}$, and
we have the following formula which is true always:
$$\begin{array}{ll}(*)\;\;\xi_r(M,e_r)=&(2\varepsilon(r)\sqrt r)^{-sign |H|}
(e_r^2-e_r^{-2})^{-2+sign |H|}\times\\
&\cases{2{r\over c}\varepsilon(c_1)\sqrt {c_1}\varepsilon(c_2)\sqrt {c_2}G,
& if $D=0$\cr
2\varepsilon(c_1)\sqrt {c_1}\varepsilon(c_2)\sqrt {c_2}d
\varepsilon(r/_{\f{cd}})\sqrt {\dfrac{r}{cd}}G,
& if $D\not=0$}\end{array}$$
for some algebraic integer   $G$

{\bf Case 1.} $(c_1,c_2) > 1$. By Lemma 4.16 in [LL2], it appears that
$$(\varepsilon(r)\sqrt r)^{-1}\varepsilon(c_1)\sqrt {c_1}\varepsilon(c_2)
\sqrt {c_2}{r\over c}=\pm\varepsilon({c_1c_2}/_{\f{c}})
\sqrt{\dfrac{c_1c_2}{c}}\varepsilon(r_{\f{c}})\sqrt {r\over c}$$
and
$$\varepsilon(c_1)\sqrt {c_1}\varepsilon(c_2)\sqrt {c_2}
\varepsilon(r/_{\f{cd}})
\sqrt{\dfrac{r}{cd}}d=\pm\varepsilon({c_1c_2}/_{\f{c}})
\sqrt{\dfrac{c_1c_2}{c}}\varepsilon(r)\sqrt {r}\varepsilon(d)
\sqrt {d}$$
Since $(c_1,c_2)>1$, we have $c_1>1,c_2>1$ and $\frac{c_1c_2}{c}> 1$.
Then by Corollary 5.4 in [LL1],
$$
\dfrac{\varepsilon(c_i)\sqrt {c_i}}
{e_r^2-e_r^{-2}}, \;i=1,2,\;\dfrac{\varepsilon(
{c_1c_2}/_{\f{c}})\sqrt{\dfrac{c_1c_2}{c}}}{e_r^2-e_r^{-2}}\;\,
\mbox{and}\;\varepsilon(b)\sqrt {b}\;\mbox{for some positive factor}\;b\;
\mbox{of}\; r$$
are all algebraic integers, so is $\xi_r(M,e_r)$ by $(*)$.
\bigskip

{\bf Case 2.} $(c_1,c_2)=1$. In this case, both $S^+$ and $S^-$ are nonempty.
\smallskip

{\bf Case 2.1.}  $D=0$. If $\dfrac{r}{c}>1$, since $\dfrac{r}{c}=\pm(\sqrt{\dfrac{r}{c}}
\varepsilon({r}/_{\f{c}}))^2$, we are done by $(*)$ and the fact
$({e_r^2-e_r^{-2}})^{-1}\sqrt{r\over c}\varepsilon(r/_{\f{c}})$
being an algebraic integer.

If $r=c$, then
$\sum_{x=1}^{r/_{\f{c}}}e_{r/_{\f{c}}}^{c(Dx^2+E_{\eta,s}n)}=1$, and
$$F_2^+ + F_2^-=2\sum^m_{s=1}(e_r^{\beta_{+},s} -e_r^{\beta_{-},s})$$
Therefore
$$\xi_r(M,e_r)=(2\varepsilon(r)\sqrt r)^{-sign |H|}
(e_r^2-e_r^{-2})^{-2+sign |H|}\varepsilon(c_1)\sqrt{c_1}\varepsilon(c_2)\sqrt{c_2}
2\sum^m_{s=1}(e_r^{\beta_{+},s} -e_r^{\beta_{-},s})$$
Since for any $a,b \in Z, e_r^a-e_r^b$ can be written as
$e_r^u(e_r^{2v}-e_r^{-2v})$ for some $u, v\in Z$, and $c=r>1$
implies that one of $c_1$ and $c_2$ must $>1$, we are done for  the case of
$H=0$. If $H\not= 0$, we have
$$\begin{array}{ll}(\varepsilon(r)\sqrt r)^{-1}
\varepsilon(c_1)\sqrt{c_1}\varepsilon(c_2)\sqrt{c_2}&=
\pm(\varepsilon(c)\sqrt c)^{-1}\varepsilon(c_1c_2)\sqrt{c_1c_2}=
\pm \varepsilon({c_1c_2}/_{\f{c}})\sqrt{\frac{c_1c_2}{c}}\\
&=\pm 1
\end{array}$$
It is done again.
\bigskip

{\bf Case 2.2}. $D\not=0$. Then by a formula in Case 1,
$$\begin{array}{ll}\xi_r(M,e_r)=&\pm(2\varepsilon(r)\sqrt r)^{-sign |H|}
(e_r^2-e_r^{-2})^{-2+sign |H|}
\varepsilon({c_1c_2}/_{\f{c}})\sqrt{\dfrac{c_1c_2}{c}}\\
&\varepsilon(r)\sqrt {r}\varepsilon(d)
\sqrt {d}\sum^m_{s=1}(e_r^{\beta_{+,s}-F_{+,s}}-e_r^{\beta_{-,s}-F_{-,s}})
\end{array}$$
Writing $e_r^a-e_r^b$  as $e_r^u(e_r^{2\nu}-e_r^{-2\nu})$,
for some $u,v \in Z$ again, we are done.

Notice that in the case $H=0,\;\nu(M)=1$, and $\xi_r(M,e_r)=2$ times an
algebraic integer. So $\Theta_r(M,e_{2r})=2^{-\nu(M)}\xi_r(M,-e_{2r})$ is an
algebraic integer.
Theorem 2 and Remark 2 are proved.
\bigskip

\section {Calculation for $X_{tref}(0)$}
\hskip\parindent
We have by the Corollary of Theorem 1 that if $3\mid\!\llap /r$, then
$$\begin{array}{ll}
\xi_r(X(3/_{\f{1}},6/_{\f{1}}, -2/_{\f{1}}), e_r)=&e_r^{-12(s^\surd (1,3)+s^
\surd (1,6))}(e_r^2-e_r^{-2})^{-2}\times\\
&\sum^r_{j=1}((e_r^{2\cdot 3'j}-e_r^{-2\cdot 3'j})
(e_r^{2\cdot 6'j}-e_r^{-2\cdot 6'j})\sum^{2'-1}_{k=0}e^{2j(2'-1-2k)}_r)
\end{array}$$
since $H=0,\;|P|=6^2$, and $s(1,-2)=0$. Expand the term $\sum^r_{j=1}$ to
$$\begin{array}{ll}
\sum^{2'-1}_{k=0}\sum^r_{j=1}(e_r^{2j(3'+6'+2'-1-2k)}
&-e_r^{2j(3'-6'+2'-1-2k)}\\&-e_r^{2j(-3'+6'+2'-1-2k)}+e_r^{2j(-3'-6'+2'-1-2k)})
\end{array}$$
We should look for the solutions of following equations for $0\leq k\leq
2'-1=\dfrac{r-1}{2}$:
$$\begin{array}{ll}
(1)\;2(3'+6'+2'-1-2k)\equiv 0&\pmod r\cr
(2)\;2(3'-6'+2'-1-2k)\equiv 0&\pmod r\cr
(3)\;2(-3'+6'+2'-1-2k)\equiv 0&\pmod r\cr
(4)\;2(-3'-6'+2'-1-2k)\equiv 0&\pmod r\cr
\end{array}
$$
Multiplying the equations by $3$, they become
$$\begin{array}{ll}
(1)\;-12k\equiv 0&\pmod r\cr
(2)\;-2-12k\equiv 0&\pmod r\cr
(3)\;-4-12k\equiv 0&\pmod r\cr
(4)\;-6-12k\equiv 0&\pmod r\cr
\end{array}$$

Thus (1) has just one solution $k=0$. For (2), $-6k=1 \pmod r$. We should
look for the solutions of $l$ satisfying $1\leq 2l+1\leq r$ and
$$\dfrac{r-(2l+1)}{2}\times (-6)\equiv 1$$
i.e. $3(2l+1)=1+nr$. Then $n$ can be only $2$, and only when $r\equiv 1 \pmod 3$, there
is a solution. For (3), we have the same conclusion as for (2). For (4),
$2k\equiv -1$ has a unique solution
$k=\frac{r-1}{2}$.
$$-12s(1,3)=\dfrac{-2}{3}\;\;\mbox{and}\;\;-12s(1,6)=\dfrac{-10}{3}$$
So
$$\xi_r(X_{tref}(0),e_r)=\cases{0, &if $r\equiv 1 \pmod 3$\cr
e_r^{-4} \dfrac{2r}{(e_r^2-e_r^{-2})^2},&if $r\equiv -1\;\;\, \pmod 3$\cr}
$$
and
$$\tau^{'}_r(X_{tref}(0))=\cases {0 &if  $r\equiv 1 \pmod 3$\cr
\sqrt{\dfrac{1}{r}}\sin{\dfrac{\pi}{r}}
\dfrac{2r}{(e_r^{2'}-e_r^{-2'})^2}e_r^{-1}=
{{-\sqrt r}\over{2\sin{\pi\over r}}}e_r^{-1},
&if $r\equiv -1 \pmod 3$\cr}$$

 \begin{center} {\bf\Large \bf   References}\end{center}
 \begin{description}

\item{} [KM1]   Kirby, R., Melvin, P.: { The 3-manifold invariants of Witten
and Reshetikhin-Turaev for $sl(2,C)$}, Invent. Math., {\bf 105} (1991) 473-545
\item{} [KM2]   Kirby, R., Melvin, P.: {Dedekind sums, $\mu$-invariants and the
signature cocycle,} Math. Ann.,{\bf 299} (1994) 231-267
\item{} [La]   Lawrence, R.: { Asymptotic expansions of Witten-
Reshetikhin-Turaev invariants for some simple 3-manifolds}, J. Math. Phys.,
{\bf 36} (1995) 6106-6129
\item{} [Li]   Li, B.H.: { Relations among Chern-Simons-Witten-Jones invariants
}, Science in China, series A, {\bf 38} (1995) 129-146
\item{} [LL1]  Li, B.H., Li, T.J.:{ Generelized Gaussian Sums and Chern-Simons-
Witten-Jones invariants of Lens spaces}, J. Knot theory and its Ramifications,
{\bf 5}  (1996) 183-224
\item{} [LL2]  Li,B.H., Li, T.J.: Kirby-Melvin's $\tau_{r}^{'}$ and Ohtsuki's $\tau$
for lens spaces, Chinese Sci. Bull.,{\bf 44} (1999)  423-426
\item{} [LL3]  Li B.H., Li, T.J.: Does Ohtsuki invariant determine full quantum
$SO(3)-$ invariants?  To appear in Letters in Math. Phys.
\item{} [LL4]  Li B.H., Li, Q.S.: Witten invariants of plumbed 3-manifolds
  Chinese Math. Ann., series A, {\bf 17} (1996) 565-572
\item{} [O]  Ohtsuki, T.: { A polynomial invariant of rational homology
3-spheres}, Inv. Math., {\bf 123} (1996) 241-257
\item{} [R1]  Rozansky, L.:Witten's invariant of rational homology spheres at
prime values of K and trivial connection contribution, Commun. Math. Phys.,
{\bf 180} (1996) 297-324.
\item{} [R2]  Rozansky, L.: { On p-adic properties of the Witten-
Reshetikhin-Turaev invariants}, math., QA/9806075, 12 Jun. (1998)
\item{} [RT]  Reshetikhin, N.Yu., Turaev, V.G.: { Invariants of 3-manifolds
via link polynomials and quantum groups}, Invent. Math., {\bf 103} (1991)
547-597
\item{} [W]  Witten, E: { Quantum field theory and the Jones polynomial},
Comm. Math.Phys., {\bf 121} (1989) 351-399
\end{description}
\vspace{.3cm}

 Author's address:\\
 Bang-He Li\\
             Institute of Systems Science\\
             Academia Sinica \\
    Beijnig 100080\\
    P. R. China \\
     Libh@iss06.iss.ac.cn \\
\end{document}